\documentclass[a4paper,11pt]{amsart}

\usepackage{a4,amsmath,amssymb,amscd,verbatim}
\newtheorem{prop}[equation]{Proposition}

\newtheorem{thm}[equation]{Theorem}
\newtheorem{cor}[equation]{Corollary}
\newtheorem{lem}[equation]{Lemma}

\theoremstyle{definition}

\newtheorem{defn}[equation]{Definition}
\newtheorem{rem}[equation]{Remark}

\newtheorem{exa}[equation]{Example}
\newtheorem{exas}[equation]{Examples}

\numberwithin{equation}{section}

\newcommand{\sands}{\mbox{$\quad\text{and}\quad$}}

\newcommand{\sts}[1]{\mbox{$\quad\text{#1}\quad$}}
\newcommand{\mstms}[1]{\mbox{$\;\;\text{#1}\;\;$}}

\newcommand{\coker}{\operatorname{coker}}

\newcommand{\ip}[1]{\mbox{$\langle{#1}\rangle$}}
\newcommand{\GL}{\mbox{${\it GL}$}}

\newcommand{\BU}{\mbox{\it BU}}

\newcommand{\bC}{\mathbb{C}}

\newcommand{\bR}{\mathbb{R}}
\newcommand{\bZ}{\mathbb{Z}}

\newcommand{\cod}{\operatorname{cod}}

\newcommand{\brpn}{\mbox{$\bR_{\scriptscriptstyle\geqslant}^n$}}

\newcommand{\brpm}{\mbox{$\bR_{\scriptscriptstyle\geqslant}^m$}}

\newcommand{\srf}[1]{\mbox{${\text{\it SR}^F}$}}

\newcommand{\daaja}{Davis and Januszkiewicz}

\newcommand{\bbox}{\mathbin{\Box}}
\newcommand{\bcs}{\mathbin{\#}}

\newcommand{\zp}{\mathcal Z_P}
\newcommand{\zq}{\mathcal Z_Q}

\def\geq{\geqslant}
\def\le{\leqslant}
\def\leq{\leqslant}

\begin{document}
\bibliographystyle{plain}
\title
{Spaces of polytopes and cobordism of quasitoric
manifolds}
\author{Victor M Buchstaber}
\address{Steklov Mathematical Institute,
Russian Academy of Sciences, Gubkina Street 8, 119991 Moscow,
Russia} \email{buchstab@mendeleevo.ru}
\author{Taras E Panov}
\address{Department of Mathematics and Mechanics, Moscow State
University, Leninskie Gory, 119992 Moscow, Russia}
\email{tpanov@mech.math.msu.su}
\author{Nigel Ray} \address{School of
Mathematics, University of Manchester, Oxford Road, Manchester
M13~9PL, England}
\email{nige@ma.man.ac.uk}

\thanks{The first and second authors were supported by the Russian
Foundation for Basic Research, grants number 04-01-00702 and
05-01-01032. The second author was supported by an EPSRC Visiting
Fellowship at the University of Manchester}



\begin{abstract}
Our aim is to bring the theory of analogous polytopes to bear on
the study of quasitoric manifolds, in the context of stably
complex manifolds with compatible torus action. By way of
application, we give an explicit construction of a quasitoric
representative for every complex cobordism class as the quotient
of a free torus action on a real quadratic complete intersection.
We suggest a systematic description for omnioriented quasitoric
manifolds in terms of combinatorial data, and explain the
relationship with non-singular projective toric varieties
(otherwise known as toric manifolds). By expressing the first and
third authors' approach to the representability of cobordism
classes in these terms, we simplify and correct two of their
original proofs concerning quotient polytopes; the first relates
to framed embeddings in the positive cone, and the second involves
modifying the operation of connected sum to take account of
orientations. Analogous polytopes provide an informative setting
for several of the details.
\end{abstract}

\noindent\emph{To Askold Khovanskii, a brilliant mathematician and
pioneer of toric geometry, on the occasion of his 60th birthday}

\smallskip

\maketitle

%
%
%
%
%
%
%
%
%

\section{Introduction}\label{intr}

The theory of analogous polytopes was initiated by Alexandrov
\cite{al37} in the 1930s, and extended more recently by Khovanskii
and Pukhlikov \cite{pu-kh92}. Our aim is to apply this theory to
the algebraic topology of torus actions.

\daaja\ \cite{da-ja91} explain how to construct a
$2n$--dimensional manifold $M$ from a characteristic pair
$(P,\lambda)$, where $P$ is a simple convex polytope of dimension
$n$, and $\lambda$ is a function with certain special properties
which assigns a subcircle of the torus $T^n$ to each facet of $P$.
By construction $M$ admits a locally standard $T^n$ action, whose
quotient space is homeomorphic to $P$. \daaja\ describe such
manifolds as {\it toric}; more recently, the term {\it quasitoric}
has been adopted, to avoid confusion with the non-singular compact
toric varieties of algebraic geometry. We follow this convention
below, and refer to such $M$ as {\it quasitoric manifolds}.

Every simple polytope $P$ is equivalent to an arrangement
$\mathcal{H}$ of $m$ closed half-spaces in an $n$--dimensional vector
space $V$, whose bounding hyperplanes meet only in general
position. The intersection of the half-spaces is assumed to be
bounded, and defines $P$. The $(n-1)$--dimensional faces $F_j$ are the
facets of $P$, where $1\le j\le m$, and general position ensures that
any face of codimension $k$ is the intersection of precisely $k$
facets. In particular, every vertex is the intersection of $n$ facets,
and lies in an open neighbourhood isomorphic to the positive cone
$\bR^n_{\geq}$. For any characteristic pair $(P,\lambda)$, it is
possible to vary $P$ within its combinatorial equivalence class
without affecting the equivariant diffeomorphism type of the
quasitoric manifold $M$.

For a fixed arrangement, we consider the vector $d_\mathcal{H}$ of
signed distances from the origin $O$ to the bounding hyperplanes
in $V$; a coordinate is positive when $O$ lies in the interior of
the corresponding half-space, and negative in the complement. We
then identify the $m$--dimensional vector space $\bR^m$ with the
space of arrangements {\it analogous to $\mathcal{H}$}. Under this
identification, $d_\mathcal{H}$ corresponds to $\mathcal{H}$
itself, and every other vector corresponds to the arrangement
obtained by the appropriate parallel displacement of half-spaces.
For small displacements, the intersections of the half-spaces are
polytopes similar to $P$. For larger displacements the
intersections may be degenerate, or empty; in either case, they
are known as {\it virtual polytopes}, and are {\it analogous to
$P$}.

In \cite{bu-ra01}, the first and third authors consider
dicharacteristic pairs $(P,\ell)$, where $\lambda$ is replaced by
a homomorphism $\ell\colon T^m\rightarrow T^n$. This has the
effect of orienting each of the subcircles $\lambda(F_j)$ of
$T^n$, and leads to the construction of an omnioriented quasitoric
manifold $M$; \cite[Theorem 3.8]{bu-ra01} claims that a canonical
stably complex structure may then be chosen for $M$. The proof,
however, has two flaws. Firstly, it fails to provide a
sufficiently detailed explanation of how a certain complexified
neighbourhood of $P$ may be framed, and secondly, it requires an
orientation of $M$ (and hence of $P$) for the stably complex
structure to be uniquely defined. The latter issue has already
been raised in \cite[\S 5.3]{bu-pa02}, but amended proofs have not
appeared. One of our aims is to show that analogous polytopes
offer a natural setting for several of the missing details.

The main application of \cite[Theorem 3.8]{bu-ra01} is as follows.

\smallskip
\noindent {\bf Theorem \ref{6.11}.}\;\cite[Theorem 6.11]{bu-ra01}
\emph{In dimensions $>2$, every complex cobordism class contains a
quasitoric manifold, necessarily connected, whose stably complex
structure is induced by an omniorientation, and is therefore
compatible with the action of the torus.}

\smallskip

This result builds upon a construction~\cite{bu-ra98} of a special
set of additive generators for the complex cobordism groups
$\varOmega_n^U$, represented by quasitoric manifolds. The proof
proceeds by considering $2n$--dimensional omnioriented quasitoric
manifolds $M_1$ and $M_2$, with quotient polytopes $P_1$ and $P_2$
respectively, and constructs a third such manifold $M$, which is
complex cobordant to the connected sum $M_1\bcs M_2$. For the
quotient polytope of $M$, the authors use the connected sum
$P_1\bcs P_2$, over which the dicharacteristics naturally extend.

In the light of our preceding observations, we must amend this
proof to incorporate the orientations of $P_1$ and $P_2$. However,
it is not always possible to form $P_1\bcs P_2$ in the oriented
sense, and simultaneously extend the dicharacteristics. Instead,
we replace $M_2$ with a complex cobordant quasitoric manifold
$M_2'$, whose quotient polytope is $I^n\bcs P_2$, where $I^n$
denotes an appropriately oriented $n$--cube. It turns out that the
resulting gain in geometrical freedom allows us to extend both the
orientations and the dicharacteristics; the result is the
omnioriented quasitoric manifold $M_1\bcs M_2'$ over the polytope
\[
P_1\bbox P_2\;=\;P_1\bcs I^n\bcs P_2,
\]
which we call the {\it box sum\/} of $P_1$ and $P_2$. We may then
complete the proof of Theorem~\ref{6.11} as described in Section
\ref{cosu} below.

In dimension $2$, $P_1\bbox P_2$ is combinatorially equivalent to the
Minkowski sum $P_1+P_2$, which is central to the theory of analogous
polytopes.

In \cite{bu-ra01}, the authors compare Theorem~\ref{6.11} with a
famous question of Hirzebruch, who asks for a description of those
complex cobordism classes which may be represented by connected
algebraic varieties. This is a difficult problem, and remains
unsolved; nevertheless, our modification to the proof of
Theorem~\ref{6.11} adds some value to the comparison, in the following
sense.

Given complex cobordism classes $[N_1]$ and $[N_2]$ of the same
dimension, suppose that $N_1$ and $N_2$ are connected. Then we may
form the connected sum $N_1\bcs N_2$ in the standard fashion, so
that it is also a connected stably complex manifold, and
represents $[N_1]+[N_2]$. If, on the other hand, $N_1$ and $N_2$
are algebraic varieties, then $N_1\bcs N_2$ is not usually
algebraic. In these circumstances we might proceed by analogy with
the quasitoric case, and look for an alternative representative
$N_2'$ such that $N_1\bcs N_2'$ is also algebraic.

We now outline the contents of each section, with additional comments
where appropriate.

In Section \ref{anpo} we recall various definitions and notation
concerning simple convex polytopes with ordered facets. We
introduce the space $\bR(P)$ of polytopes analogous to a fixed
example $P$, and consider a linear map $\chi_P\colon V\to\bR(P)$,
defined on the ambient space $V$ of $P$. Under $\chi_P$, a point
$y\in V$ is mapped to that point of $\bR(P)$ which represents the
polytope congruent to $P$ obtained by shifting the origin to~$y$.
We then interpret the projection from $\bR(P)$ to the cokernel of
$\chi_P$ as mapping a polytope $P'\in\bR(P)$ to the vector of
distances from a distinguished vertex $v_\star\in P'$ to its $m-n$
opposite facets. This allows us to describe the projection
explicitly, as a map $C\colon\bR(P)\to\bR^{m-n}$.

In Section \ref{toma} we summarise the construction of a
quasitoric manifold $M$ over a polytope $P$ with $m$ facets
\cite{da-ja91}. In their work, Davis and Januszkiewicz use an
auxiliary $T^m$-space $\zp$, whose quotient by the kernel of a
dicharacteristic homomorphism they identify with $M$. It has
become evident that the spaces $\zp$ are of great independent
interest in toric topology, and they are now known as
\emph{moment-angle complexes}~\cite{bu-pa02}. They arise in
homotopy theory as homotopy colimits \cite{p-r-v03}, in symplectic
topology as level surfaces for the moment maps of Hamiltonian
torus actions~\cite{pano??}, and in the theory of arrangements as
complements of coordinate subspace arrangements~\cite[Chapter
8]{bu-pa02}. Using the matrix of the projection $C$, we describe
$\zp$ as a complete intersection of real quadratic hypersurfaces.

In Section \ref{scorfr} we amend the definition of omniorientation so
as to include an orientation of $M$, and recall the stably complex
structure which results. In so doing, we frame $\zp$ equivariantly in
$\bR^{2m}$ and consider the quotient framing on $P$ as a submanifold
of the positive cone in the space of analogous polytopes.

We review the construction of connected sum for omnioriented
quasitoric manifolds in Section \ref{cosu}, by encoding the
additional orientations as signs attached to the fixed points. We
then explain how to correct \cite[Theorem 3.8]{bu-ra01}, and
recover Theorem \ref{6.11}. Combining the latter with our
quadratic description of $\zp$ yields the following additional
result on representability.

\smallskip

\noindent {\bf Theorem \ref{repquad}. }\emph{Every complex cobordism
class may be represented by the quotient of a free torus action on a
real quadratic complete intersection.}

\smallskip

The importance of the real quadratic viewpoint has been emphasised in
recent work of Bosio and Meersseman~\cite{bo-me??}, who consider a
specific class of complete intersections of real quadrics in $\mathbb
C^m$, called \emph{links}. They show that all links (taking products
with a circle in odd-dimensional cases) can be endowed with the
structure of a non-K\"ahler complex manifold, generalising the class
of Hopf and Calabi-Eckmann manifolds. It is clear from both the
results of~\cite{bo-me??} and our own Section~\ref{toma} that the
class of links coincides with the class of moment-angle complexes
$\zp$ arising from simple polytopes. This fact provides further
connections between toric topology and complex geometry, in which
calculations of cohomology rings of moment-angle complexes carried out
in~\cite{bu-pa02} feature prominently.

Finally, in Section \ref{excore}, we discuss the realisation of
$4$--dimensional complex cobordism classes by omnioriented quasitoric
manifolds, and comment on comparable situations in higher dimensions.

Throughout our work we adopt the combinatorial convention that
$[n]$ denotes the set of integers $\{1,\dots,n\}$, for any natural
number $n$. Occasionally, it is convenient to interpret $[0]$ as
the empty set. We write $2^{[n]}$ for the Boolean algebra of
subsets of $[n]$, ordered by inclusion or its reverse as
necessary.

We are pleased to acknowledge the input of several colleagues in
preparing this article. In particular, Konstantin Feldman stimulated
our development of the box sum by observing that the equivariant
connected sum of~\cite{bu-ra01} and~\cite{da-ja91} cannot be used to
produce a quasitoric representative for the cobordism class $2[\mathbb
C P^2]$. Mikiya Masuda supplied valuable assistance in our
understanding of Example~\ref{nspv}, and Peter McMullen offered
helpful guidance on simple polytopes.  Tony Bahri and Neil Strickland
provided additional criticisms of \cite{bu-ra01}, which alerted its
authors to the need for clarification.
%
%
%
%
%
%
%
%
%

\section{Analogous polytopes}\label{anpo}

We work in a real vector space $V$ of dimension $n$, equipped with
a euclidean inner product $\ip{\;\;,\;\;}$ and an orthonormal
basis $e_1$, \dots, $e_n$. An {\it ordered arrangement}
$\mathcal{H}$ of {\it closed half-spaces\/} in $V$ is a collection
of subsets
\begin{equation}\label{hasp}
H_i=\{x\in V:\ip{a_i,x}+b_i\geq0\}\sts{for}1\leq i\leq m,
\end{equation}
where $a_i$ lies in $V$ and $b_i$ is a real scalar. Unless stated
otherwise, we assume that $\mathcal{H}$ has cardinality $m\geq n$, and
that $a_i$ has unit length for every $1\leq i\leq m$. We consider
$H_i$ to be a smooth manifold, whose boundary $\partial H_i$ is its
{\it bounding hyperplane}
\begin{equation}\label{bohy}
Y_i=\{x\in V:\ip{a_i,x}+b_i=0\}\sts{for}1\leq i\leq m,
\end{equation}
with inward pointing normal vector $a_i$

When the intersection $\cap_iH_i$ is bounded, it forms a {\it convex
polytope} $P$; otherwise, it is a {\it polyhedron}. We assume that $P$
has maximal dimension $n$ and that $\mathcal{H}$ is irredundant, in
the sense that no $H_i$ may be deleted without enlarging $P$. In these
circumstances, $\mathcal{H}$ and $P$ are interchangeable. We may also
specify $P$ by a matrix inequality $A_Px+b_P\geq0$, where $A_P$ is the
$m\times n$ matrix of row vectors $a_i$, and $b_P$ is the column
vector of scalars $b_i$ in $\bR^m$. If we permute the half-spaces
\eqref{hasp} by an element of the symmetric group $\varSigma_m$, we
recover $P$ by applying the same permutation to the rows of $A_P$ and
the coordinates of $b_P$.

\begin{exas}\label{conesplx}
The {\it standard $n$--simplex $\varDelta(n)$} is the polytope defined
by the half-spaces
\begin{equation}\label{splx}
H_i=\begin{cases}
\{x:\ip{e_i,x}\geq 0\}&\text{for $1\leq i\leq n$}\\
\{x:\ip{a_{n+1},x}+1\geq0\}&\text{for $i=n+1$}
\end{cases}
\end{equation}
in $\bR^n$, where $a_{n+1}=(-1,\dots,-1)$; its vertices are the points
$0$, $e_1$, \dots, $e_n$. The {\it positive cone\/} $\brpn$ is the
polyhedron obtained by deleting $H_{n+1}$ from \eqref{splx}; it has a
single vertex $0$, and contains all vectors with non-negative
coordinates.
\end{exas}

\begin{exa}\label{ncubepoly}
The {\it standard $n$--cube $I^n$} is the polytope defined by the
half-spaces
\begin{equation}\label{ncube}
H_i=\begin{cases}
\{x:\ip{e_i,x}\geq 0\}&\text{for $1\leq i\leq n$}\\
\{x:-\ip{e_i,x}+1\geq 0\}&\text{for $n+1\leq i\leq 2n$}
\end{cases}
\end{equation}
in $\bR^n$; its vertices are the binary sequences
$(\delta_1,\dots,\delta_n)$, where $\delta_j=0$ or $1$ for $1\leq
j\leq n$.
\end{exa}

A {\it supporting hyperplane\/} is characterised by the property
that $P$ lies within one of its two associated half-spaces. A {\it
proper face\/} of $P$ is defined by its intersection with any
supporting hyperplane, and forms a convex polytope of lower
dimension. We regard $P$ as an $n$--dimensional face of itself;
the faces of dimension $0$, $1$, and $n-1$ are known as {\it
vertices}, {\it edges}, and {\it facets} respectively. There is
one facet $F_i=P\cap Y_i$ for every bounding hyperplane
\eqref{bohy}, so the facets correspond bijectively to the
half-spaces \eqref{hasp}. We deem a vertex $v$ and facet $F_i$ to
be {\it opposite\/} whenever $v$ lies in the interior of $H_i$. If
the bounding hyperplanes are in general position, then every
vertex of $P$ is the intersection of exactly $n$ facets, and has
$m-n$ opposite half-spaces. In these circumstances, $P$ is {\it
simple}.

From this point on, we deal only with simple polytopes, and
reserve the notation $q=q(P)$ and $m=m(P)$ for the number of
vertices and facets respectively. Every face of codimension $k$
may be written uniquely as
\begin{equation}\label{facelabs}
F_I\;=\;F_{i_1}\cap\dots\cap F_{i_k}
\end{equation}
for some subset $I=\{i_1,\dots,i_k\}$ of $[m]$, and the $F_I$ may
then be ordered lexicographically for each $1\leq k\leq n$.

The faces are elements of the {\it face poset} $\mathcal{L}(P)$,
ordered by reverse inclusion. We use the subscripts $I$ to interpret
$\mathcal{L}(P)$ as a subposet of the Boolean algebra $2^{[m]}$,
ordered by inclusion; so $\mathcal{L}(P)$ is {\it ranked\/} \cite[p
99]{stan86} by the codimension function $\cod(F_I)=|I|$. It has unique
minimal element $F_\varnothing=P$, and its maximal elements are the
vertices. It fails to be a lattice only because we usually omit the
empty face, which would otherwise form a unique maximal element.

Two polytopes are {\it combinatorially equivalent\/} whenever
their face posets are isomorphic.
A combinatorial equivalence class of polytopes is known as a {\it
combinatorial polytope}, and most of our constructions are defined
on such classes. Nevertheless, it often helps to keep a
representative polytope in mind, rather than the underlying poset.
Natural examples of combinatorial polytopes include the {\it
vertex figures} $P_v$, which are formed by intersecting $P$ with
any closed half-space whose interior contains a single vertex $v$.
Because $P$ is simple, $P_v$ is an $n$--simplex for any $v$.

By permuting the facets of $P$ if necessary, we may assume that
the intersection $F_1\cap\dots\cap F_n$ is a vertex $v_\star$. In
this case, we describe $P$ as {\it finely ordered}, and refer to
$v_\star$ as the {\it initial vertex\/}; it is the first vertex of
$P$ with respect to the ordering implied by \eqref{facelabs}. For
computational purposes it is often convenient to locate the
initial vertex of $P$ at the origin, and use the normal vectors
$a_i$ as an orthonormal basis for $V$. This may be achieved by
applying an affine transformation to $P$, which preserves its
combinatorial type.

Given a second finely ordered polytope $P'$ in $\bR^{n'}$, we may
list the facets of the product $P\times P'$ as
\begin{equation}\label{fiorpr}
F_1\times P',\;\dots,\; F_m\times P',\; P\times F_1',\;\dots,\;
P\times F'_{m'},
\end{equation}
where $F_i$ and $F'_j$ range over the facets of $P$ and $P'$
respectively. Then $P\times P'$ is finely ordered by shifting the
block of facets $P\times F_1'$, \dots, $P\times F'_{n'}$ to the left,
until it occupies positions $n+1$, \dots, $n+n'$ in
\eqref{fiorpr}. The initial vertex is $(v_\star,v'_\star)$. Of course,
this procedure yields different results for $P\times P'$ and $P'\times
P$.

For a fixed arrangement $\mathcal{H}$, we consider the vector
$d_\mathcal{H}\in\bR^m$, whose $i\/$th coordinate is the signed
distance from $Y_i$ to the origin $O$ in $V$, for $1\leq i\leq m$.
The sign is positive when $O$ lies in the interior of $H_i$, and
negative in the exterior. So long as we maintain our convention
that the normal vectors $a_i$ have unit length, $d_\mathcal{H}$
coincides with $b_P$; otherwise, the distances have to be scaled
accordingly. Every vector $d_\mathcal{H}+h$ in $\bR^m$ may then be
identified with an {\it analogous arrangement\/} of half-spaces,
defined by translating each $H_i$ by $h_i$, for $1\leq i\leq m$.
Some such arrangements determine convex polytopes $P(h)$, and
others, dubbed {\it virtual polytopes}, do not. In either case,
they are deemed to be {\it analogous to $P$}. We note that $P(h)$
is given by
\begin{equation}\label{defph}
\{x\in V:A_Px+b_P+h\geq0\},
\end{equation}
and is combinatorially equivalent to $P$ when $h$ is small. In
particular, we have that $P(0)=P$.

\begin{exas}\label{basepolys}
The zero vector $0\in\bR^m$ is identified with the central
arrangement $\mathcal{H}_0$, whose bounding hyperplanes contain
the origin; the corresponding polytope $P(-b_P)=\{0\}$ is virtual.
The basis vector $e_i\in\bR^m$ is identified with the arrangement
obtained from $\mathcal{H}_0$ by translating $H_i$; the
corresponding polytope $P(-b_P+e_i)=P_i$ may be virtual, or a
simplex.
\end{exas}

\begin{exa}\label{congpolys}
Any $y\in V$ defines a vector $A_Py\in\bR^m$. Then $A_Py+b_P$ is
identified with the arrangement obtained by translating
$\mathcal{H}$ by $-y$; the corresponding polytope $P(A_Py)$ is the
translate $P-y$, and is congruent to $P$. As $y$ varies, we obtain
an $n$--parameter family of analogous polytopes, each being
congruent to $P$.
\end{exa}

The {\it Minkowski sum\/} of subsets $P$, $Q\subseteq V$ is given
by
\[
P+Q=\{x+y:x\in P,\,y\in Q\}\;\subseteq\;V.
\]
If $P$ and $Q$ are convex polytopes, so is $P+Q$; moreover, when
$P$ is analogous to $Q$, so is $P+Q$. Under the identification of
$b_P+h$ with $P(h)$, vector addition corresponds to Minkowski sum,
and scalar multiplication to rescaling. In this context, we denote
the $m$--dimensional vector space of polytopes analogous to $P$ by
$\bR(P)$, and consider the identification as an isomorphism
\begin{equation}\label{khid}
k\colon\bR^m\longrightarrow \bR(P),\sts{where} k(b_P+h)=P(h).
\end{equation}

We may interpret the matrix $A_P$ as a linear transformation
$V\rightarrow \bR^m$. Since the points of $P$ are specified by the
constraint $A_Px+b_P\geq0$, the intersection of the affine
subspace $A_P(V)+b_P$ with the positive cone $\brpm$ is a copy of
$P$ in $\bR^m$. In other words, the formula $i_P(x)=A_Px+b_P$
defines an affine injection
\begin{equation}\label{defip}
i_P\colon V\longrightarrow\bR^m,
\end{equation}
which embeds $P$ as a submanifold of the positive cone. Since $i_P$
maps the half-space $H_i$ to the half-space $\{y:y_i\geq 0\}$, it
embeds each codimension--$k$ face of $P$ in a codimension--$k$ face of
$\brpm$.

The composition $\chi_P=k\circ i_P$ restricts to an affine injection
$P\rightarrow \bR(P)$, and Example \eqref{congpolys} identifies
$\chi_P(y)$ as the polytope congruent to $P$, obtained by translating
the origin to $y$, for all $y$ in $P$. Of course, $\chi_P(P)$ is a
submanifold of the positive cone
$\bR(P)_{\scriptscriptstyle\geqslant}$, and facial codimensions are
preserved as before.

When $P$ is finely ordered, the half spaces $H_1+h_1$, \dots,
$H_n+h_n$ determine the initial vertex $v_\star(h)$ of $P(h)$ for
any shift vector $h$. For every $1\leq i\leq m$, we write $d_i(h)$
for the signed distance between $v_\star(h)$ and the supporting
hyperplane $Y_i+h_i$; in other words,
\[
  d_i(h)=\ip{a_i,v_\star(h)}+b_i+h_i\mstms{for all}1\leq i\leq m,
\]
and $d_i(h)=0$ for $1\leq i\leq n$. We define a linear
transformation $C\colon\bR^m\rightarrow\bR^{m-n}$ by the formula
\begin{equation}\label{codef}
C(b_P+h)\;=\;(d_{n+1}(h),\ldots,d_{m}(h)).
\end{equation}
Using \eqref{khid}, we may interpret $C$ as a transformation
$\bR(P)\rightarrow \bR^{m-n}$, which acts by
$P(h)\mapsto(d_{n+1}(h),\ldots,d_{m}(h))$. Clearly $C$ is
epimorphic.

\begin{prop}\label{nuzfr}
As a transformation $V\rightarrow \bR^{m-n}$, the composition
$C\cdot A_P$ is zero.
\end{prop}
\begin{proof}
The $d_i(h)$ are metric invariants of the polytope $P(h)$, so $C$
takes identical values on congruent polytopes. In particular, it
is constant on the translates $P-y$ for all values $y\in V$, and
therefore on the affine plane $A_P(V)+b_P$. So $C(A_P(V))=0$, as
required.
\end{proof}
Proposition \ref{nuzfr} determines a short exact sequence of the
form
\begin{equation}\label{ses}
  0\longrightarrow V\stackrel{A_P}{\longrightarrow}\bR^{m}
  \stackrel{C}{\longrightarrow}\bR^{m-n}\longrightarrow 0,
\end{equation}
or equivalently, a choice of basis for $\coker A_P$.

In order to construct a matrix $(c_{i,j})$ for $C$, it is convenient
to use the orthonormal basis $a_1$, \dots, $a_n$, as described
above. Then the basis polytopes $P_j$ of \eqref{basepolys} satisfy
\[
  d_i(P_j)=
  \begin{cases}
    -a_{i,j}&\text{if $1\le j\le n$}\\
    \delta_{i,j}&\text{if $n+1\le j\le m$}
  \end{cases}
\]
for all $n+1\leq i\leq m$, giving
\begin{equation}\label{coma}
(c_{i,j})\;=\;
\begin{pmatrix}
  -a_{n+1,1} & \ldots & -a_{n+1,n} & 1     & 0      & \ldots & 0\\
  -a_{n+2,1} & \ldots & -a_{n+2,n} &      0 & 1     & \ldots & 0\\
  \vdots    & \ddots & \vdots    & \vdots &\vdots&\ddots&\vdots\\
  -a_{m,1}   & \ldots & -a_{m,n}   & 0 & 0  & \ldots & 1
\end{pmatrix}.
\end{equation}

A permutation of the facets produces an alternative basis for
$\coker A_P$, and the corresponding matrix is obtained by
reordering the columns of $C$. Any other $(m-n)\times m$ matrix $C$
for which $CA_P=0$ also provides a basis for $\coker A_P$, so long
as it has full rank; it necessarily satisfies the following
property.
\begin{lem}\label{minor}
Let $C'$ be the $(m-n)\times(m-k)$ matrix obtained from $C$ by
deleting columns $c_{j_1}$, \dots, $c_{j_k}$, for some $1\leq
k\leq n$; if the intersection $F_{j_1}\cap\dots\cap F_{j_k}$ is a
face of $P$ of codimension $k$, then $C'$ has rank $m-n$.
\end{lem}
\begin{proof}
Let $\iota\colon\bR^{m-k}\rightarrow\bR^m$ be the inclusion of the
subspace
\[
\{x:x_{j_1}=\dots=x_{j_k}=0\}
\]
and $\kappa\colon\bR^m\rightarrow\bR^k$ the associated quotient
map. Then $C'$ is the matrix of the composition $C\cdot\iota$, and
the $k\times n$ matrix $A'$ of the composition $\kappa\cdot A_P$
consists of the rows $a_{j_1}$, \dots, $a_{j_k}$ of $A_P$. The
data implies that $A'$ has rank $k$, and therefore that
$\kappa\cdot A_P$ is an epimorphism; so $C\cdot\iota$ is also an
epimorphism, and its matrix has rank $m-n$.
\end{proof}
%
%
%
%
%
%
%
%
%

\section{Quasitoric manifolds}\label{toma}

In this section we include a summary of \daaja's construction of
quasitoric manifolds $M$ over a simple polytope~$P$. Throughout,
we use the methods and notation of \cite{bu-ra01}, and refer
readers to \cite[Chapter 6]{bu-pa02}) for futher details. We also
assume that $P$ is finely ordered; so $M$ has a distinguished
fixed point, near which we insist that $T^n$ act as standard. To
illustrate these additional requirements we revisit several
standard examples.

We denote the $i\/$th coordinate subcircle of the standard
$m$--torus $T^m$ by $T_i$ for every $1\leq i\leq m$. Given any
subset $I\subseteq [m]$, we define the subgroup $T_I$ by
\[
\prod_{i\,\in\, I}T_i\;<\;T^m,
\]
so that $T_\varnothing$ is the trivial subgroup $\{1\}$. Every
point $p$ of $P$ lies in the interior of a unique face $F_{I_p}$,
where $I_p$ is given by $\{i:p\in F_i\}$; we abbreviate $F_{I_p}$
and $T_{I_p}$ to $F(p)$ and $T(p)$ respectively. If $p$ is a
vertex, for example, then $T(p)$ has dimension $n$ (the maximum
possible), and if $p$ is an interior point of $P$, then $T(p)$ is
trivial.

We define the {\it moment-angle complex\/} $\zp$ as the identification
space
\begin{equation}\label{defzp}
  T^m\times P/{\sim}\;,
\end{equation}
where $(t_1,p)\sim(t_2,p)$ if and only if $t_1^{-1}t_2\in T(p)$.  So
$\zp$ admits a canonical left $T^m$-action whose isotropy subgroups
are precisely the subgroups $T(p)$. Construction \eqref{defzp} may
equally well be applied to the positive cone $\brpm$, in which case
the result is the complex vector space $\bC^m$. Since the embedding
$i_P$ of \eqref{defip} respects facial codimensions, there is a
pullback diagram
\begin{equation}\label{cdiz}
\begin{CD}
  \zp @>i_Z>>\bC^m\\
  @V\varrho_P VV\hspace{-0.2em} @VV\varrho V @.\\
  P @>i_P>> \bR^m_\geq
\end{CD}
\end{equation}
of identification spaces. Here $\varrho(z_1,\ldots,z_m)$ is given by
$(|z_1|^2,\ldots,|z_m|^2)$, the vertical maps are projections onto the
quotients by the $T^m$-actions, and $i_Z$ is a $T^m$-equivariant
embedding. It is sometimes convenient to rewrite $\bC^m$ as
$\bR^{2m}$, in which case we substitute $q_j+ir_j$ for the $j\/$th
coordinate $z_j$, and let $T_j$ act by rotation.

Then Proposition \ref{nuzfr} and Diagram \eqref{cdiz} imply that
$i_Z$ embeds $\zp$ in $\bR^{2m}$ as the space of solutions of the
$m-n$ quadratic equations
\begin{equation}\label{zpquad}
  \sum_{k=1}^mc_{j,k}\left(q_k^2+r_k^2-b_k\right)=0,\sts{for}
  1\leq j\leq m-n.
\end{equation}
In Lemma \ref{izfr}, we will confirm that $\zp$ is a frameable
submanifold of $\bR^{2m}$ of dimension $(m+n)$, and therefore smooth.

In order to construct quasitoric manifolds over $P$, we need one
further set of data. This consists of a homomorphism $\ell\colon
T^m\rightarrow T^n$, satisfying \daaja's independence condition,
namely
\begin{equation}\label{dajastar}
F_I\;\;\text{is a face of codimension $k$}\quad
\Longrightarrow\quad \text{$\ell\,$ is monic on $T_I$}.
\end{equation}
Any such $\ell$ is called a {\it dicharacteristic\/} in
\cite{bu-ra01}. Condition \eqref{dajastar} ensures that the kernel
$K(\ell)$ of $\ell$ is isomorphic to an $(m-n)$--dimensional
subtorus of $T^m$, and features in a short exact sequence
\begin{equation}\label{ses2}
  1\longrightarrow K(\ell)\longrightarrow T^m
  \stackrel{\ell}{\longrightarrow}T^n\longrightarrow 1.
\end{equation}
Wherever possible we abbreviate $K(\ell)$ to $K$.

We write the subcircle $\ell(T_i)<T^n$ as $T(F_i)$ for any $1\leq
i\leq m$, and the subgroup $\ell(T_I)$ as $T(F_I)$ for any face
$F_I$. For each point $p$ in $P$ we let $S(p)$ denote the subgroup
$T(F(p))$; it is, of course, $\ell(T(p))$. For example, $S(w)=T^n$
for any vertex $w$, and $S(p)=\{1\}$ for any point $p$ in the
interior of $P$.

When applied to the initial vertex $v_\star$, \eqref{dajastar}
ensures that the restriction of $\ell$ to $T_1\times\dots\times
T_n$ is an isomorphism. So we may use the circles $T(F_1)$, \dots,
$T(F_n)$ to define a basis for the Lie algebra of $T^n$, and
represent the homomorphism induced by $\ell$ by an $n\times m$
integral matrix of the form
\begin{equation}\label{lamat}
\varLambda\;=\;\begin{pmatrix}
  1&0&\ldots&0&\lambda_{1,n+1}&\ldots&\lambda_{1,m}\\
  0&1&\ldots&0&\lambda_{2,n+1}&\ldots&\lambda_{2,m}\\
  \vdots&\vdots&\ddots&\vdots&\vdots&\ddots&\vdots\\
  0&0&\ldots&1&\lambda_{n,n+1}&\ldots&\lambda_{n,m}
\end{pmatrix}.
\end{equation}
It is often convenient to partition $\varLambda$ as
$\left(I_n\;|\;\varLambda_\star\right)$, so that $\varLambda_\star$ is
an $n\times(m-n)$-matrix. Given any other vertex
$F_{j_1}\cap\dots\cap F_{j_n}$, \eqref{dajastar} implies that the
corresponding columns $\lambda_{j_1}$, \dots, $\lambda_{j_n}$ form
a basis for $\bZ^n$, and have determinant $\pm1$. We refer to
\eqref{lamat} as the {\it refined form\/}, and call
$\varLambda_\star$ the {\it refined submatrix\/} of $\ell$.

Since $K$ meets every isotropy subgroup $T(p)$ of the $T^m$-action
trivially, it acts freely on $\zp$, and the base of the resulting
principal $K$-bundle $\pi_\ell\colon\zp\rightarrow M$ is a smooth
$2n$--dimensional manifold. By construction, $M$ may be expressed as
the identification space
\begin{equation}\label{deftoma}
T^n\times P/{\approx}
\end{equation}
where $(s_1,p)\approx(s_2,p)$ if and only if $s_1^{-1}s_2\in
S(p)$. Furthermore, $M$ admits a canonical $T^n$-action $\alpha$,
which is locally isomorphic to the standard action on $\bC^n$, and
has quotient map $\pi\colon M\rightarrow P$. Note that
$\pi\cdot\pi_\ell=\varrho_P$ as maps $\zp\rightarrow P$. The fixed
points of $\alpha$ project to the vertices of $P$, so they are
also ordered, and we refer to $\pi^{-1}(v_\star)$ as the {\it
initial fixed point} $x_\star$. Then \eqref{deftoma} identifies a
neighbourhood of $x_\star$ with $\bC^n$, on which $\alpha$ is
standard; its representation at other fixed points may be read off
from the corresponding columns of $\varLambda$.

The quadruple $(M,\alpha,\pi,P)$ is an example of a {\it
quasitoric manifold\/}, as defined by \daaja. Any manifold with a
similarly well-behaved torus action over $P$ is
$\theta$-equivariantly homeomorphic to one of the form
\eqref{deftoma}, see~\cite[Prop.~1.8]{da-ja91}. In this sense, $M$
is typical, and we follow the lead of \cite{bu-ra01} by assuming
that every quasitoric manifold is presented in the form
\eqref{deftoma}.

Additional structure on $M$ is associated to the {\it facial
submanifold} $M_i$, defined as the inverse images of the facet
$F_i$ under $\pi$, for $1\leq i\leq m$. It is clear that every
$M_i$ has codimension $2$, and that its isotropy subgroup is
$T(F_i)<T^n$. The quotient map
\begin{equation}\label{defrhoh}
\zp\times_K\bC_i\longrightarrow M
\end{equation}
defines a canonical complex line-bundle $\rho_i$, whose
restriction to $M_i$ is isomorphic to the normal bundle $\nu_i$ of
its embedding in $M$.

The submanifolds $M_i$ are mutually transverse, and we write
\begin{equation}\label{facemans}
M_I\;=\;M_{i_1}\cap\dots\cap M_{i_k}
\end{equation}
for any non-empty intersection, using $I$ as in \eqref{facelabs}.
So $M_I$ is the inverse image of the codimension--$k$ face $F_I$
under $\pi$. Of course $M_I$ has codimension $2k$, and its
isotropy subgroup is $T(F_I)$. The restriction of
$\rho_I=\rho_{i_1}\oplus\dots\oplus\rho_{i_d}$ to $M_I$ is
isomorphic to the normal bundle $\nu_I$ of its embedding in $M$,
for any face $F_I$.

As explained in \cite {da-ja91}, the bundles $\rho_i$ play an
important part in understanding the integral cohomology ring of
$M$. If $u_i$ denotes the first Chern class $c_1(\rho_i)$ in
$H^2(M)$, then $H^*(M)$ is generated by $u_1$, \dots, $u_m$,
modulo two sets of relations. The first are linear, and arise from
the refined form \eqref{lamat} of the dicharacteristic; the second
are monomial, and arise from the Stanley-Reisner ideal of $P$. The
former may be read off from the refined submatrix as
\begin{equation}\label{ui}
  u_i\;=\;-\lambda_{i,n+1}u_{n+1}-\ldots-\lambda_{i,m}u_m
  \sts{for}1\leq i\leq n,
\end{equation}
and show that $u_{n+1}$, \dots, $u_m$ suffice to generate $H^*(M)$
multiplicatively.

In work such as \cite{bu-pa02}, \cite{bu-ra01}, and
\cite{ci-ra05}, fine orderings are not considered, so the matrices
representing $\ell$ rarely appear in refined form. In order to
rectify this situation systematically, we may begin by choosing an
initial vertex. Then we shuffle the facets of $P$ (and therefore
the columns of the representing matrix) until $F_{[n]}$ is
$v_\star$, and premultiply the resulting matrix by the unique
element of $\GL(n;\bZ)$ that transforms the first $n$ columns into
$I_n$. To illustrate the effects of this procedure, we now revisit
three important examples.

\begin{exa}\label{kostyaex1}
The $n$--simplex is finely ordered by \eqref{splx}, and has
initial vertex the origin. Then $i_P$ embeds $\varDelta(n)$ in
$\bR^{n+1}$ by $i_P(x)=(x_1,\dots,x_n,1-\sum_{i=1}^nx_i)$, and
$\zp$ is the unit sphere $S^{2n+1}\subset\bC^{n+1}$. The refined
submatrix is the column vector $(-1,\dots,-1)$ in $\bR^n$, so the
kernel of the dicharacteristic is the diagonal subcircle
\[
T_\delta\;=\;\{(t,t,\dots,t)\}\;<\;T^{n+1}.
\]
It follows that $M$ is the complex projective space $\bC P^n$.
Then $(t_1,\dots,t_n)\in T^n$ acts on the point with homogeneous
coordinates $[z_1,\dots,z_{n+1}]$ as multiplication by
$(t_1,\dots,t_n,1)$, and the initial fixed point is
$[0,\dots,0,1]$. Every facial bundle is isomorphic to
$\overline{\eta}$, where $\eta$ is the Hopf line bundle. The
cohomology ring of $M$ is generated by elements $u_1$, \dots,
$u_{n+1}$ in $H^2(M)$, and relations \eqref{ui} give
$u_1=\cdots=u_{n+1}$; the Stanley-Reisner relations reduce to
$u_1^{n+1}=0$.
\end{exa}

\begin{exa}\label{botttow1}
The $n$--cube is finely ordered by \eqref{ncube}, and has initial
vertex the origin. Then $i_P$ embeds $I^n$ in $\bR^{2n}$ by
$i_P(x)= (x_1,\dots,x_n,1-x_1,\dots,1-x_n)$, and $\zp$ is the
product of unit $3$--spheres $|z_k|^2+|z_{n+k}|^2=1$ in
$\bC^{2n}$, where $1\leq k\leq n$. The refined submatrix is
\[
D\;=\; {\scriptsize
\begin{pmatrix}
-1&0&\dots&0&0&0&\dots&0&0\\
d(1,2)&-1&\dots&0&0&0&\dots&0&0\\
&&\vdots&&&&\vdots&&\\
d(1,j)&d(2,j)&\dots&d(j-1,j)&-1&0&\dots&0&0\\
&&\vdots&&&&\vdots&&\\
d(1,n)&d(2,n)&\dots&d(j-1,n)&d(j,n)&d(j+1,n)& \dots&d(n-1,n)&-1
\end{pmatrix}
}
\]
for any set of $n(n-1)/2$ integers $d(i,j)$, where $1\leq i<j\leq
n$; so the kernel of the dicharacteristic is the $n$--torus
\[
\{(t_1,t_1^{-d(1,2)}t_2,\dots,t_1^{-d(1,n)}t_2^{-d(2,n)}\dots
t_{n-1}^{-d(n-1,n)}t_n,t_1,t_2,\dots,t_n)\}\;<\;T^{2n}.
\]
It follows that $M$ is the $n$th stage $Q_n$ of the {\it Bott
tower\/} defined in \cite{ci-ra05} and \cite{gr-ka94}, albeit with
permuted coordinates. Then $(t_1,\dots,t_n)\in T^n$ acts on the
equivalence class $[z_1,\dots,z_{2n}]$ as multiplication by
$(t_1,\dots,t_n,1,\dots,1)$, and the initial fixed point is
$[0,\dots,0,1,\dots,1]$. The facial bundles are the $\rho_i$ of
\cite{ci-ra05}, suitably reordered. The cohomology ring of $M$ is
generated by $u_1$, \dots, $u_{2n}$ in $H^2(M)$, and relations
\eqref{ui} give
\[
u_j\;=\;-d(1,j)u_{n+1}-\dots-d(j-1,j)u_{n+j-1}+u_{n+j}
\sts{for}1\leq j\leq n.
\]
The Stanley-Reisner relations take the form $u_ju_{n+j}=0$ for all
$j$.

If the defining integers satisfy
\[
d(i,j)=
\begin{cases}
1&\text{for $i=j-1$}\\
0&\text{otherwise}
\end{cases}
\]
for every $2\leq j\leq n$, then $Q_n$ becomes the {\it bounded
flag manifold\/} $B_n$ of \cite{bu-ra98}. If $d(i,j)=0$ for all
$i$, $j$, then $Q_n$ is the $n$--fold product $(S^2)^n$.
\end{exa}

\begin{exa}\label{brs1}
For any pair of integers $r<s$, the facets of
$R^{r+s-1}=I^r\times\varDelta(s-1)$ are finely ordered by
combining \eqref{fiorpr} with Examples \ref{kostyaex1} and
\ref{botttow1}. The initial vertex is the origin.  Then $i_P$
embeds $R^{r+s-1}$ in $\bR^{2r+s}$ by
\[
i_P(x)= \Big(x_1,\dots,x_{r+s-1},1-x_1,\dots,
1-x_r,1-{\textstyle\sum_{i=r+1}^{r+s-1}}x_i\Big),
\]
and $\zp$ is a product $S^3\times\dots\times S^3\times S^{2s-1}$
of $r+1$ unit spheres in $\bC^{2r+s}$. The refined submatrix
(which is $(r+s-1)\times(r+1)$) is
\[
E\;=\;\begin{pmatrix}
J_r&0\\
J_r&-1\\
0_{r,s}&-1
\end{pmatrix},
\]
where $J_r$ is the $r\times r$ matrix whose only non-zero elements are
$-1\/$s on the diagonal and $1\/$s on the subdiagonal, $0_{r,s}$ is
the $(s-r-1)\times r$ zero matrix, and $0$, $-1$ denote column vectors
of the appropriate length. So the kernel of the dicharacteristic is
the $(r+1)$--dimensional subtorus
\[
\{(t_1,t_1^{-1}t_2,\dots,t_{r-1}^{-1}t_r,ut_1,ut_1^{-1}t_2,\dots,
ut_{r-1}^{-1}t_r,u,\dots,u,t_1,t_2,\dots,t_r,u)\}
\]
of $T^{2r+s}$. It follows that $M$ is the $\bC P^{s-1}$-bundle
$B_{r,s}$ over $B_r$ as defined in \cite{bu-ra98}, albeit with
permuted coordinates. Then $(t_1,\dots,t_{r+s-1})\in T^{r+s-1}$ acts
on the equivalence class $[z_1,\dots,z_{2r+s}]$ as multiplication by
$(t_1,\dots,t_{r+s-1},1,\dots,1)$, and the initial fixed point is
$[0,\dots,0,1,\dots,1]$. The facial bundles are the $\rho_i$ of
\cite{ci-ra05}, suitably reordered. The cohomology ring of $B_{r,s}$
is generated by elements $u_1$, \dots, $u_{2r+s}$ in $H^2(M)$, and it
is helpful to write $u_{r+s+i-1}=x_i$ for $1\leq i\leq r$, and
$u_{2r+s}=y$. Then relations \eqref{ui} give
\begin{multline*}
u_1=x_1,\quad u_{r+1}=x_1+y,\quad u_{2r+1}=\dots=u_{r+s-1}=y,\\
\text{and}\quad u_i=x_i-x_{i-1},\quad u_{r+i}=x_i-x_{i-1}+y\quad
\text{for}\quad 2\leq i\leq r.
\end{multline*}
The Stanley-Reisner relations are of the form $u_ix_i=0$ for
$1\leq i\leq r$, and $u_{r+1}u_{r+2}\dots u_{r+s-1}y=0$. Thus
$H^*(B_{r,s})$ is isomorphic to
\[
\bZ[x_1,\dots,x_r,y]\big/J,\sts{where} J=\big( x_i^2-x_ix_{i-1},
\,{\textstyle\sum_{j=0}^r}x_r^{j}y^{s-j}\big).
\]
The same construction works for $r=s$, but $\varLambda_\star$ is more
complicated.
\end{exa}

%
%
%
%
%
%
%
%
%

\section{Stably complex structures, orientations, and
framings}\label{scorfr}

On a smooth manifold $N$ of dimension $d$, a {\it stably complex
structure\/} is an equivalence class of real $2k$--plane bundle
isomorphisms $\tau(N)\oplus\bR^{2k-d}\cong\zeta$. Here $\zeta$
denotes a fixed $\GL(k,\bC)$--bundle, $\bR^{2k-d}$ denotes the
trivial $(2k-d)$--dimensional bundle with fibre $\bR^{2k-d}$ and
$k$ is suitably large. Two such isomorphisms are equivalent when
they agree up to stabilisation; or, alternatively, when the
corresponding lifts to $\BU$ of the classifying map of the stable
tangent bundle of $N$ are homotopic through lifts. Note that
$\bR^{2k-d}$ is canonically oriented (and even framed) by choosing
the standard basis, which therefore determines an orientation
for~$N$.

Now assume that $N$ has an $l$--dimensional torus action
$\alpha\colon T^l\times N\to N$. A stably complex structure on $N$
is \emph{$T^l$--invariant} whenever the composition
\[
\begin{CD}
  \zeta @>\cong >> \tau(N)\oplus\bR^{2k-d}
  @>\;d\alpha(t,\:\cdot)\oplus 1\;>>
  \tau(N)\oplus\bR^{2k-d} @>\cong>> \zeta
\end{CD}
\]
is an isomorphism of complex bundles for every $t\in T^l$. In this
section we show that every quasitoric manifold admits an invariant
stably complex structure and identify the geometric data required
to induce these structures.

According to \cite{bu-ra01}, an {\it omniorientation\/} of a
quasitoric manifold $M$ consists of a choice of orientation for
each normal bundle $\nu_i$. This coincides with a choice of
complex structure for each $\rho_i$, and is therefore equivalent
to a dicharacteristic $\ell$. In \cite{bu-pa02}, a choice of
orientation for $M$ is also assumed, since none is implied by
$\ell$. We adopt this convention henceforth, and refer to the
constituent data as the dicharacteristic and orientation {\it
associated to\/} the omniorientation. The orientation corresponds
to a fundamental class $\mu_M$ in the integral homology group
$H_{2n}(M)$.

An interior point of the quotient polytope $P$ admits an open
neighborhood $U$, whose inverse image under the projection $\pi$
is canonically diffeomorphic to $T^n\times U$ as a subspace of
$M$. Since $T^n$ is oriented by the choice of basis leading to the
refined form~\eqref{lamat} of the matrix of $\ell$, orientations
of $M$ correspond bijectively to orientations of $P$. Every pair
$(P,\varLambda)$ therefore determines a $2n$--dimensional
omnioriented quasitoric manifold, where $P$ is the combinatorial
type of an oriented finely ordered $n$--dimensional simple
polytope, and $\varLambda$ is a matrix of the form \eqref{lamat}.

\begin{defn}\label{data}
We refer to the pair $(P,\varLambda_\star)$ as the {\it combinatorial
data\/} underlying the omnioriented manifold $M$.
\end{defn}

We may specify the orientation of $P$ on a representative polytope
in $\bR^n$, or by an equivalence class of orderings of the $n$
edges incident on $v_\star$ in $\mathcal{L}(P)$. The latter is
independent of the fine ordering on $P$ (although they may, of
course, agree). When it is important to emphasise that the facial
submanifolds of $M$ are ordered, and that $\alpha$ is standard at
$x_\star$, we also describe $M$ as {\it refined}.

In order to explain the stably complex structure induced on $M$,
it is convenient to study the embedding $i_Z$ of \eqref{cdiz} in
more detail.
\begin{lem}\label{izfr}
The embedding $i_Z\colon\zp\rightarrow\bR^{2m}$ is
$T^m$--equivariently framed by any choice of matrix $C=(c_{i,j})$
for the transformation~\eqref{codef}.
\end{lem}
\begin{proof}
We describe $i_Z$ by the $m-n$ quadratic equations \eqref{zpquad}
over $P\subset \brpm$. At each point
$(q_1,r_1,\dots,q_m,r_m)\in\zp$, the $m-n$ associated gradient
vectors are given by
\begin{equation}\label{grve}
2\left(c_{j,1}q_1,\,c_{j,1}r_1,\,\dots,\,c_{j,m}q_m,\,c_{j,m}r_m\right)
\sts{for}1\leq j\leq m-n,
\end{equation}
and so form the rows of the $(m-n)\times 2m$ matrix $2CR$, where
\[
R\;=\;
\begin{pmatrix}
  q_1&r_1& \ldots & 0&0 \\
  \vdots &\vdots & \ddots & \vdots & \vdots\\
  0&0  & \ldots &q_m&r_m
\end{pmatrix}
\]
is $m\times 2m$. By definition of $i_P$, the set of integers
$j_1$,\dots, $j_k$ with the property that
$q_{j_1}=r_{j_1}=\dots=q_{j_k}=r_{j_k}=0$ at some point $z\in\zp$
corresponds to an intersection $F_{j_1}\cap\dots\cap F_{j_k}$ of
facets forming a face of $P$ of codimension $k$. Lemma \ref{minor}
then applies to show that the matrix obtained by deleting the
columns $c_{j_1}$, \dots, $c_{j_k}$ of $C$ has rank $m-n$. It
follows that $2CR$ has rank $m-n$, and therefore that the gradient
vectors \eqref{grve} are linearly independent at $z$, and so frame
$i_Z$.

Furthermore, each of the gradient vectors frames the corresponding
quadratic hypersurface in $\bR^{2m}$, and is $T^m$--invariant.
\end{proof}

\begin{rem}
Lemma \ref{izfr} provides an alternative to \cite[Proposition
3.4]{bu-ra01}, where insufficient detail is given for readers to
complete the proof.
\end{rem}

It is particularly illuminating to describe the framing of $i_Z$
in terms of analogous polytopes, as follows.

Factoring out by the action of $T^m$ yields a framing of the
embedding $i_P$, and therefore of $P$ in $\brpm$; moreover, on
each face of $P$, the framing lies in the ambient face of $\brpm$.
Under the identification \eqref{khid}, the framing vectors may be
represented by $m-n$ independent $1$--parameter families of
polytopes analogous to $P$. These families are made explicit by
applying the differential $d\varrho_P$ to the rows of the matrix
$2CR$. At the point $(q_1,r_1,\dots,q_m,r_m)$ in $\zp$, the matrix
of $d\varrho_P$ is given by $2R$, so the framing vectors are the
rows of the $(m-n)\times m$ matrix $4CRR^t$. When $C$ takes the
form \eqref{coma}, we may take the $j\/$th framing vector to be
\[
f_j\;=\;
(-a_{n+j,1}y_1,\dots,-a_{n+j,n}y_n,0,\dots,0,y_{n+j},0,\dots,0)
\]
at $y=i_P(x)$, for $1\leq j\leq m-n$. Applying \eqref{khid}, we
conclude that the corresponding $1$--parameter family of polytopes
$P(f_j,t)$ (for $-1\leq t\leq 1$) is obtained from $P$ by:
retaining the origin at $x$, rescaling $H_k$ by $-a_{n+j,k}t$ for
$1\leq k\leq n$, fixing every facet opposite the initial vertex
except $H_{n+j}$, and rescaling the latter by $t$.

It is possible to reverse this procedure, and begin by framing
$i_P$. The corresponding $T^m$--equivariant framing of $i_Z$ is then
recovered by applying the contruction \eqref{defzp}. Since $P$ is
contractible, all framings of $i_P$ are equivalent, and their lifts to
$i_Z$ are equivariantly equivalent. In particular, the equivalence
class of the framings described in Lemma \ref{izfr} does not depend on
the choice of fine ordering on $P$.

The smoothness of $M$ is assured by Lemma \ref{izfr}, and we now
return to its tangent bundle $\tau(M)$. Our analysis reduces to a
special case of Szczarba's proof of \cite[Theorem (1.1)]{szcz64}, and
supercedes that given in~\cite[Theorem (3.8)]{bu-ra01} which ignores
the orientation of $M$.
\begin{prop}\label{rbund}
Any omnioriented quasitoric manifold admits a canonical stably
complex structure, which is invariant under the $T^n$--action.
\end{prop}
\begin{proof}
There is a $T^m$--equivariant decomposition
\[
  \tau(\zp)\oplus\nu(i_Z)\cong\zp\times\bC^m,
\]
obtained by restricting the tangent bundle $\tau(\bC^m)$ to $\zp$.
Factoring out by the kernel of $\ell\colon T^m\to T^n$ yields
\begin{equation}\label{tauzp}
  \tau(M)\oplus\left(\xi/K\right)\oplus
  \left(\nu(i_Z)/K\right)\cong
  \zp\times_{K}\bC^m,
\end{equation}
where $\xi$ denotes the $(m-n)$--plane bundle of tangents along
the fibres of $\pi_\ell$. The right-hand side of \eqref{tauzp} is
isomorphic to $\bigoplus_{i=1}^m\rho_i$ as $\GL(m,\bC)$--bundles.

Szczarba \cite[Corollary 6.2]{szcz64} identifies $\xi/K$ with the
adjoint bundle of $\pi_\ell$, which is trivial because $K$ is
abelian; and $\nu(i_Z)/K$ is trivial by Lemma \ref{izfr}. So
\eqref{tauzp} reduces to an isomorphism
\begin{equation}\label{tansumrho}
\tau(M)\oplus\bR^{2(m-n)}\cong\rho_1\oplus\ldots\oplus\rho_m,
\end{equation}
although different choices of trivialisations may lead to
different isomorphisms. Since $M$ is connected and
$\GL(2(m-n),\bR)$ has two connected components, such isomorphisms
are equivalent when and only when the induced orientations agree
on $\bR^{2(m-n)}$. We choose the orientation which is compatible
with those for $\tau(M)$ and $\rho_1\oplus\ldots\oplus\rho_m$, as
given by the omniorientation.

The induced structure is invariant under the action of $T^n$,
because $i_Z$ is $T^m$--equivariant.
\end{proof}

The stably complex structures represented by the two choices of
orientation differ by sign. The underlying smooth structure is
also $T^n$--invariant, and is identical to that inferred from Lemma
\ref{izfr}.

The proof of Proposition \ref{rbund} allows us to evaluate the
tangential Chern classes of the canonical stably complex
structure.
\begin{cor}\label{mcherns}
In $H^{2i}(M)$, the Chern class $c_i(\tau)$ is given by the $i$th
symmetric polynomial in the variables $u_1$, \dots, $u_m$, for
$1\leq i\leq n$.
\end{cor}
\begin{proof}
By \eqref{tansumrho}, the total Chern class of $\tau$ is
$c(\tau)=\prod_{i=1}^m(1+u_i)$ in $H^*(M)$.
\end{proof}

\daaja's quasitoric manifolds are inspired by the non-singular
projective toric varieties of algebraic geometry. Every such $X$ is
determined by the normal fan of an integral simple polytope
$Q\subset\bR^n$, whose vertices lie in the lattice $\bZ^n$. We may
assume that the origin is a distinguished vertex, that its incident
facets lie in the respective coordinate hyperplanes, and that the
remaining facets $F_{n+1}$, \dots, $F_m$ are ordered. Since $X$ is
equipped with a canonical complex structure it is also an omnioriented
quasitoric manifold, so we study this example in more detail before
moving on.

According to Batyrev \cite[\S8.2]{bu-pa02}, $X$ may be identified with
the \emph{geometric quotient} of the \emph{coordinate subspace
complement}
\[
  U(Q)\;=\;\bC^m\setminus\bigcup\,\{z:z_{i_1}=\ldots=z_{i_k}=0
  \text{ if }F_{i_1}\cap\ldots\cap F_{i_k}=\varnothing\text{ in }Q\}
\]
by the complexified group
\[
  K_\bC=\ker(\ell_\bC\colon(\bC^*)^m\to(\bC^*)^n).
\]
By definition, there is a canonical embedding
$j\colon\zq\stackrel{\subset}{\longrightarrow}U(Q)$ of the compact
subset $\zq$, which induces an algebraic isomorphism $\zq/K\rightarrow
U(Q)/K_\bC$. In other words, $\zq$ is the {\it Kempf-Ness set\/} for
the action of the algebraic torus $K_\bC$ on the quasiaffine variety
$U(Q)$; see~\cite[Theorem~3.4]{pano??} for further details. The
integral version
\[
  0\longrightarrow V_{\bZ}\stackrel{A_Q}{\longrightarrow}\bZ^{m}
  \stackrel{C}{\longrightarrow}\bZ^{m-n}\longrightarrow 0
\]
of the short exact sequence \eqref{ses} is the sequence of weight
lattices for \eqref{ses2}.

The tangent bundle of $U(Q)$ is trivial, and admits a
$\GL(m-n,\bC)$--subbundle $\xi_\bC$ of tangents along the fibre of
the quotient map $U(Q)\to X$. Applying \cite{szcz64} once more, we
deduce that the complex structure on $X$ is compatible with the
corresponding isomorphism
\begin{equation}\label{taux}
  \tau(X)\oplus(\xi_\bC/K_\bC)\;\cong\;U(Q)\times_{K_\bC}\bC^m
\end{equation}
of quotient $\GL(m,\bC)$--bundles, where $\xi_\bC/K_\bC$ is trivial
because $K_\bC$ is abelian.

\begin{exa}\label{nspv}
For any non-singular projective toric variety $X$, we let $P$ be the
oriented combinatorial type of $Q$, and the columns of $\varLambda$ be
the primitive integral inward pointing normal vectors to $F_1$, \dots,
$F_m$ respectively.  So $\varLambda=A_Q^t$ in the notation of Section
\ref{anpo} (although the row vectors $a_i$ of $A_Q$ do not necessarily
have unit length). We identify the stably complex structure associated
to the combinatorial data $(P,\varLambda_\star)$ by comparing
isomorphisms \eqref{tansumrho} and \eqref{taux} as follows.

On fibres, $j$ restricts to the inclusion of the subgroup $K<K_\bC$. So
there is an isomorphism $\xi\oplus\nu(i_Z)\cong j^*\xi_\bC$ over $\zq$,
whose quotient is an isomorphism
\[
  \xi/K\oplus\nu(i_Z)/K\;\cong\;\xi_\bC/K_\bC
\]
of $\GL(2(m-n),\bR)$--bundles over $X$. Similarly, there is a quotient
isomorphism
\[
  \zp\times_K\bC^m\;\cong\;U(Q)\times_{K_\bC}\bC^m
\]
of $\GL(m,\bC)$--bundles, which identifies the right hand sides of
\eqref{tansumrho} and \eqref{taux}. Moreover, the framing chosen for
$\xi/K\oplus \nu(i_Z)/K$ in \eqref{tansumrho} induces the same
orientation as that of the complex structure on $\xi_\bC/K_\bC$,
because both are compatible with the natural orientation of $X$. So
the stably complex structure associated to $(P,\varLambda_\star)$ agrees
with that induced by the complex varietal structure on $X$.
\end{exa}

\begin{exa}\label{kostyaex2}
The standard basis for $\bR^n$ defines an orientation of
$\varDelta(n)$; combining this with Example \ref{kostyaex1} yields the
combinatorial data $\left(\varDelta(n),-1\right)$ for $\bC P^n$. The
data actually arises from the normal fan of $\varDelta(n)$, so Example
\ref{nspv} applies, and the corresponding omniorientation agrees with
that induced by the complex structure on $\bC P^n$. The
omniorientation may be altered by conjugating the $j$th facial bundle
for any $1\leq j\leq n+1$, which has the effect of negating the $j$th
column of $\varLambda$.  For $j\leq n$, restoring the dicharacteristic
to refined form involves replacing the refined submatrix by the column
vector $\epsilon=(\epsilon_1,\dots,\epsilon_n)$, where $\epsilon_i=-1$
for $i\neq j$, and $\epsilon_j=1$. This procedure may be extended to
any subset $J\subseteq [n+1]$. We write the result as $P_\epsilon$, to
emphasise the omniorientation; the resulting stably complex structure
may be described by an isomorphism
\[
\tau(P_\epsilon)\oplus\bC\;\cong\;|J|\,\eta
\oplus\left(n+1-|J|\right)\overline{\eta}.
\]
\end{exa}

\begin{exa}\label{botttow2}
The standard basis for $\bR^n$ defines an orientation of $I^n$;
combining this with Example \ref{botttow1} yields the combinatorial
data $\left(I^n,D\right)$ for $Q_n$. The data arises from the normal
fan of a polytope combinatorially equivalent to $I^n$, so Example
\ref{nspv} applies, and the corresponding omniorientation agrees with
that induced by the complex structure on $Q_n$. The omniorientation
may be altered by conjugating the $j$th facial bundle for any $1\leq
j\leq 2n$, which has the effect of negating the $j$th column of
$\varLambda$. For $j\leq n$, restoring the dicharacteristic to refined
form involves negating the $j$th row of $D$. This procedure may be
extended to any subset $J\subseteq [2n]$, although many of the
resulting stably complex structures coincide, and several bound
\cite{ci-ra05}.

The bounding cases are no less natural to topologists than the
projective algebraic varieties, and play an important role in complex
cobordism theory \cite{ray86}. For example, when $Q_n$ is the
$n$--fold product $S=(S^2)^n$, the combinatorial data $(I^n,I_n)$
corresponds to the bounding structure given by
$\eta\oplus\overline{\eta}$ on each cartesian factor.
\end{exa}

\begin{exa}\label{brs2}
The orientations of Examples \ref{kostyaex2} and \ref{botttow2}
describe an orientation for $R^{r+s-1}$; combining this with Example
\ref{brs1} yields the combinatorial data $\left(R,E\right)$ for
$B_{r,s}$. The data arises from the normal fan of a polytope
combinatorially equivalent to $R$, so Example \ref{nspv} applies, and
the corresponding omniorientation agrees with that induced by the
complex structure on $B_{r,s}$. This stably complex structure ensures
that certain linear combinations of the $B_{r,s}$ form multiplicative
generators for the complex cobordism ring $\varOmega_*^U$
\cite{bu-ra98}.
\end{exa}

%
%
%
%
%
%
%
%
%

\section{Connected Sums}\label{cosu}

In this section we review the construction of the connected sum for
omnioriented quasitoric manifolds $M'$ and $M''$, as was sketched
in~\cite[1.11]{da-ja91} and realised in~\cite{bu-ra01}. However, the
orientations demanded by Proposition \ref{rbund} were omitted in both
descriptions, and we deal with them here in terms of signs associated
to the vertices of $P$.

We denote the dicharacteristics associated to the omniorientations
of $M'$ and $M''$ by $\ell'$ and $\ell''$, with refined
submatrices $\varLambda'_\star$ and $\varLambda''_\star$ respectively;
and assume that the associated orientations are given by
orientations of the polytopes $P'$ and $P''$ . In addition, we let
$P'$ and $P''$ be finely ordered by $o'$ and $o''$, with initial
vertices $v'_\star$ and $v''_\star$ respectively.

The {\it connected sum} $P'\bcs_{v'_\star,v''_\star} P''$ may be
described informally as follows.  First construct the polytope
$Q'$ by deleting the interior of the vertex figure $P'_{v'_\star}$
from $P'$; so $Q'$ has one new facet $\varDelta(v'_\star)$ (which
is an $(n-1)$--simplex), whose incident facets are ordered by
$o'$. Then construct the polytope $Q''$ from $P''$ by the same
procedure. Finally, glue $Q'$ to $Q''$ by identifying
$\varDelta(v'_\star)$ with $\varDelta(v''_\star)$, in such a way
that the $j\/$th facet of $Q'$ combines with the $j\/$th facet of
$Q''$ to give a single new facet for each $1\leq j\leq n$. The
gluing is carried out by applying appropriate projective
transformations to $Q'$ and $Q''$. Precise details are given in
\cite[\S6]{bu-ra01}.

The combinatorial type of the connected sum may be changed, for
example, by choosing alternative fine orderings on $P'$ and $P''$.
So long as the choices are clear, or their effect on the result is
irrelevant, we use the abbreviation $P'\bcs P''$. The face lattice
$\mathfrak{L}_F(P'\bcs P'')$ is obtained from
$\mathfrak{L}_F(P')\cup\mathfrak{L}_F(P'')$ by identifying the
$j\/$th facets of each, for $1\leq j\leq n$. In particular,
\begin{equation}\label{csstats}
q(P'\bcs P'')=q(P')+q(P'')-2\sands m(P'\bcs P'')=m(P')+m(P'')-n.
\end{equation}

By definition, the {\it connected sum}
$M'\bcs_{x'_\star,x''_\star}M''$ is the quasitoric manifold
constructed over $P'\bcs P''$ using the dicharacteristic
$\ell_\#\colon T^{m'+m''-n}\to T^n$ associated to the matrix
\begin{equation}\label{cmcs}
\varLambda_\#\;=\;\begin{pmatrix}
1&0&\ldots&0&\lambda'_{1,n+1}&\ldots&\lambda'_{1,m'}&\lambda''_{1,n+1}&
\ldots&\lambda''_{1,m''}\\
  0&1&\ldots&0&\lambda'_{2,n+1}&\ldots&\lambda'_{2,m'}&\lambda''_{2,n+1}&
\ldots&\lambda''_{2,m''}\\
  \vdots&\vdots&\ddots&\vdots&\vdots&\ddots&\vdots&\vdots&\ddots&\vdots\\
  0&0&\ldots&1&\lambda'_{n,n+1}&\ldots&\lambda'_{n,m'}&\lambda''_{n,n+1}&
\ldots&\lambda''_{n,m''}
\end{pmatrix}.
\end{equation}
Of course $\varLambda_\#$ is no longer in refined form, since the
first $n$ facets of $P'\bcs P''$ have empty intersection.
Nevertheless, we may finely order $P'\bcs P''$ by choosing the {\it
second\/} vertex of $P'$ as initial vertex, and applying the
procedure described immediately before Example \ref{kostyaex1}.

By construction, $M'\bcs M''$ is equivariantly diffeomorphic to the
equivariant connected sum of $M'$ and $M''$ at their initial fixed
points. If $M'$ and $M''$ are omnioriented, the only possible
obstruction to defining a compatible omniorientation of $M'\bcs
M''$ involves the associated orientations. We deal with this issue
in Proposition \ref{cssign} below.

We write $p'\colon M'\bcs M''\to M'$ and $p''\colon M'\bcs M''\to
M''$ for the maps collapsing the connected sum onto its
constituent manifolds.

We recall from \cite{pano99} that an omniorientation attaches a
{\it sign} $\sigma(w)$ to every vertex $w$ of the quotient
polytope $P$ (or, equivalently, to every fixed point of $M$). By
definition, $\sigma(w)=\pm 1$ measures the difference between the
orientations induced on the tangent space at $w$ by the
dicharacteristic and the fundamental class $\mu_M$ respectively.
When $w$ is $F_{i_1}\cap\dots\cap F_{i_n}$, the former is given by
the Chern class
$c_n\left(\rho_{i_1}\oplus\ldots\oplus\rho_{i_n}\right)$. So we
have that
\[
  \sigma(w)\;=\;\left\langle u_{i_1}\cdots u_{i_n},\mu_M\right\rangle.
\]

\begin{prop}\label{cssign}
The connected sum $M'\bcs_{x'_\star,x''_\star}M''$ admits an
orientation compatible with those of $M'$ and $M''$ if and only if
$\sigma(v'_\star)=-\sigma(v''_\star)$.
\end{prop}
\begin{proof}
The facets of $P'\bcs P''$ give rise to complex line bundles
$\xi_i$, $\xi'_j$ and $\xi''_k$ over $M'\bcs M''$, corresponding
to the columns of \eqref{cmcs}. We denote their first Chern
classes by
\[
c_1(\xi_i)=w_i,\quad c_1(\xi'_j)=w'_j,\sands c_1(\xi''_k)=w''_k
\]
in $H^2(M'\bcs M'')$, for
\[
1\le i\le n,\quad n+1\leq j\leq m', \sands n+1\leq k\leq m''
\]
respectively. The relations \eqref{ui} become
\[
  w_i=-\lambda'_{i,n+1}w'_{n+1}-\ldots-\lambda'_{i,m'}w'_{m'}
  -\lambda''_{i,n+1}w''_{n+1}-\ldots-\lambda''_{i,m''}w''_{m''},
\]
which imply that
\[
  w_i=p'^{\,*}u'_i+p''^{\,*}u''_i\sts{for}1\leq i\leq n.
\]
Since the first $n$ facets of $P'\bcs P''$ do not define a vertex,
it follows that $w_1\cdots w_n=0$ in $H^{2n}(M'\bcs M'')$, and
\[
  (p'^{\,*}u'_1+p''^{\,*}u''_1)\cdots
  (p'^{\,*}u'_n+p''^{\,*}u''_n)
  =p'^{\,*}(u'_1\cdots u'_n)+p''^{\,*}(u''_1\cdots u''_n)=0.
\]
For any choice of fundamental class in $H_{2n}(M'\bcs M'')$, we
deduce that
\[
  \left\langle u'_1\cdots u'_n\;,\;p'_*\mu_{M'\bcs M''}\right\rangle +
  \left\langle u''_1\cdots u''_n\;,\;p''_*\mu_{M'\bcs
  M''}\right\rangle = 0.
\]
But the corresponding orientation of $M'\bcs M''$ is compatible
with those of $M'$ and $M''$ if and only if $p'_*\mu_{M'\bcs
M''}=\mu_{M'}$ and $p''_*\mu_{M'\bcs M''}=\mu_{M''}$; that is, if
and only if
\[
  \sigma(v'_\star)+\sigma(v''_\star)=0,
\]
as required.
\end{proof}
\begin{cor}\label{scscs}
Let $M'$ and $M''$ be omnioriented quasitoric manifolds over
finely ordered polytopes $P'$ and $P''$ respectively, with
$\sigma(v'_\star)=-\sigma(v''_\star)$; then the stably complex
structure induced on $M'\bcs_{x'_\star,x''_\star}M''$ by
Proposition {\rm \ref{rbund}} and Proposition {\rm \ref{cssign}}
is equivalent to the connected sum of those induced on $M'$ and
$M''$. Moreover, the associated complex cobordism classes satisfy
\[
[M'\bcs M'']=[M']+[M''].
\]
\end{cor}
\begin{proof}
The stably complex structures on $M'$ and $M''$ combine to give an
isomorphism
\begin{equation}\label{torcs}
\begin{split}
  \tau(M'\bcs M'')\oplus\bR^{2(m'+m''-n)}\;\cong\;
  \xi_1&\oplus\ldots\oplus\xi_n\oplus\xi'_{n+1}\oplus
    \ldots\oplus\xi'_{m'}\\
  &\oplus\xi''_{n+1}\oplus\ldots\oplus\xi''_{m''}.
\end{split}
\end{equation}
As explained in \cite[Theorem 6.9]{bu-ra01}, the isomorphism
\eqref{torcs} belongs to one of the two equivalence classes
specified by Proposition \ref{rbund} over $M'\bcs M''$. The choice
of orientation is then provided by Proposition \ref{cssign}.

The equation of cobordism classes follows immediately, because the
connected sum is cobordant to the disjoint union.
\end{proof}

Proposition \ref{cssign} implies that we cannot always form the
connected sum of two omniorented quasitoric manifolds. If the sign
of every vertex of $P$ is positive, for example, then it is
impossible to construct $M\bcs M$ directly; we illustrate this
situation in Example~\ref{nspv2} below.

Corollary \ref{scscs} confirms that the complex cobordism class
$[M'\bcs M'']$ is independent of the fine orderings $o'$ and
$o''$, and therefore of the initial vertices.

\begin{exa}\label{nspv2}
For any non-singular projective toric variety, it follows from
Example \ref{nspv} that the dicharacteristic and orientation both
arise from the complex structure on $X$. So they are compatible,
and every vertex of $P$ has sign $+1$.
\end{exa}

\begin{exa}\label{prodtwos}
Example \ref{botttow2} exhibits an omniorientation of $S$, defined
by the combinatorial data $(I^n,I_n)$, which induces a bounding
stably complex structure. The signs of the vertices of $I^n$ (in
the notation of Example \ref{ncubepoly}) are given by
\[
\sigma(\delta_1,\dots,\delta_n)
=(-1)^{\delta_1}\dots(-1)^{\delta_n}.
\]
So adjacent vertices have opposite sign, and both occur with
frequency $n$.
\end{exa}

We are now in a position to prove Lemma \ref{difsi}, which
emphasises an important principle; however unsuitable a quasitoric
manifold $M$ may be for the formation of connected sums, a good
alternative representative always exists within the complex
cobordism class $[M]$.
\begin{lem}\label{difsi}
Let $M$ be an omnioriented quasitoric manifold of dimension $>2$
over a finely ordered polytope $P$; then there exists an
omnioriented $M'$ over a finely ordered polytope $P'$, such that
$[M']=[M]$ and $P'$ has at least two vertices of opposite sign.
\end{lem}
\begin{proof}
Suppose that $v_\star$ is the initial vertex of $P$. Let $S$ be
the omnioriented product of $2$--spheres of Example
\ref{prodtwos}, with initial vertex $w_\star$.

If $\sigma(v_\star)=-1$, define $M'$ to be
$S\bcs_{w_\star,v_\star}M$ over $P'=I^n\bcs_{w_\star,v_\star}P$
Then $[M']=[M]$, because $S$ bounds; moreover, adjacent pairs of
non-initial vertices of $I^n$ have opposites signs, which survive
under the formation of $P'$, as sought. If $\sigma(v_\star)=+1$,
we make the same construction using the opposite orientation of
$I^n$ (and therefore of $S$). Since $-S$ also bounds, the same
conclusions hold. In either case, $P'$ may be finely ordered as
decribed above; its initial vertex corresponds to $(0,\dots,0,1)$
in $I^n$.
\end{proof}

We may now complete the proof of our amended \cite[Theorem
6.11]{bu-ra01}.
\begin{thm}\label{6.11}
In dimensions $>2$, every complex cobordism class contains a
quasitoric manifold, necessarily connected, whose stably complex
structure is induced by an omniorientation, and is therefore
compatible with the action of the torus.
\end{thm}
\begin{proof}
  Following \cite{bu-ra01}, we consider cobordism classes $[M_1]$ and
  $[M_2]$ in $\varOmega_n^U$, represented by omnioriented quasitoric
  manifolds over quotient polytopes $P_1$ and $P_2$ respectively. It
  then suffices to construct a third such manifold $M$ such that
  $[M]=[M_1]+[M_2]$, because a set of quasitoric additive generators
  for $\varOmega_n^U$ is given by \cite{bu-ra98} for all $n>0$.

  Firstly, we follow Lemma \ref{difsi} and replace $M_2$ by $M_2'$
  over $P_2'=I^n\bcs P_2$. Then we finely order $P_2'$ so as to ensure
  that its initial vertex has opposite sign to that of $P_1$, thereby
  guaranteeing the construction of $M_1\bcs M_2'$ over $P_1\bcs P_2'$.
  The resulting omniorientation defines the required cobordism class,
  by Corollary \ref{scscs} and Lemma \ref{difsi}.
\end{proof}

We refer to the polytope $P_1\bcs I^n\bcs P_2$ of Theorem
\ref{6.11} as the {\it box sum} $P_1\bbox P_2$ of $P_1$ and $P_2$,
because of the intermediate cube. The fact that we have replaced
$P_1\bcs P_2$ by $P_1\bbox P_2$ in the proof of Theorem \ref{6.11}
does not affect the following observation of \cite{bu-ra01}: for
any complex cobordism class, the quotient polytope of a
representing quasitoric manifold may be chosen to be a connected
sum of products of simplices.

Combining Theorem~\ref{6.11} with the details of Lemma \ref{izfr}
and the quadratic description \eqref{zpquad} of $\zp$ leads to the
following interesting conclusion.
\begin{thm}\label{repquad}
Every complex cobordism class may be represented by the quotient
of a free torus action on a real quadratic complete intersection.
\end{thm}

One further deduction from Theorem \ref{6.11} is the result
of~\cite{ray86}, that every complex cobordism class contains a
representative whose stable tangent bundle is a sum of line bundles.

%
%
%
%
%
%
%
%
%

\section{Examples and concluding remarks}\label{excore}

We were taught the importance of adding an orientation to the
original definition of omniorientation by certain $4$--dimensional
examples of Feldman \cite{fe02}. In this section we describe and
develop his examples (noting that $4$ is the smallest dimension to
which Proposition~\ref{cssign} is relevant). They lead to our
concluding remarks concerning higher dimensions.

We shall use a result of~\cite{pano99}, which identifies the top Chern
number of any $2n$--dimensional omnioriented quasitoric manifold as
\begin{equation}\label{taraseq}
c_n(M)\;=\;\sum_w\sigma(w).
\end{equation}
For any quotient polytope $P$, it is also convenient to refine the
notation of \eqref{csstats} by writing
\[
q(P)=q_+(M)+q_-(M),
\]
where $q_\pm(M)$ denotes the number of vertices with sign $\pm 1$
respectively. These numbers are preserved by any $\theta$--equivariant
diffeomorphism which respects omniorientations.

When $n=2$, the complex cobordism class $[\bC P^2]$ of the
standard complex structure of Example~\ref{kostyaex1} is an
additive generator of the cobordism group
$\varOmega_4^U\cong\bZ^2$, with $c_2(\bC P^2)=3$ and $q_-(\bC
P^2)=0$. Each of the other three omniorientations of Example
\ref{kostyaex2} represents the class $[\bC P^2]-4[\bC P^1]^2$
(which is an independent additive generator), and
$q_-(P_\epsilon)$ is given by the number of negative entries in
the relevant $\epsilon$; in other words, it is $1$, $1$, or $2$.

The question then arises of representing $2[\bC P^2]$ by an
omnioriented quasitoric manifold $M$. We cannot expect to use $\bC
P^2\bcs\bC P^2$ for $M$, because no vertices of sign $-1$ are
available in $\varDelta(2)$, as required by Proposition
\ref{cssign}. Moreover, $M$ must satisfy $c_2(M)=6$, by
additivity, so the quotient polytope $P$ has $6$ or more vertices;
as observed by Feldman, it follows that $P$ cannot be
$\varDelta(2)\bcs\varDelta(2)$, which is a square! So we proceed
by appealing to Lemma \ref{difsi}, and replace the second copy of
$\bC P^2$ by the omnioriented quasitoric manifold $(-S)\bcs\bC
P^2$ over $P'=I^2\bcs\varDelta(2)$. Of course $(-S)\bcs\bC P^2$ is
cobordant to $\bC P^2$, and $P'$ is a pentagon. These observations
lead naturally to our second example.

\begin{exa}\label{cptboxcpt}
The omnioriented quasitoric manifold $\bC P^2\bcs(-S)\bcs\bC P^2$
represents $2[\bC P^2]$, and lies over the box sum
$\varDelta(2)\bbox\varDelta(2)$, which is a hexagon. Figure~$1$
illustrates the procedure diagramatically, in terms of
dicharacteristics and orientations. Every vertex of the hexagon
has sign~$1$.
\begin{figure}[h]
\begin{picture}(120,35)
  \put(0,5){\line(0,1){25}}
  \put(0,5){\line(5,3){21}}
  \put(0,30){\line(5,-3){21}}
  \put(40,5){\line(-1,1){12.5}}
  \put(40,5){\line(1,1){12.5}}
  \put(40,30){\line(-1,-1){12.5}}
  \put(40,30){\line(1,-1){12.5}}
  \put(80,5){\line(0,1){25}}
  \put(80,5){\line(-5,3){21}}
  \put(80,30){\line(-5,-3){21}}
  \put(95,10){\line(3,-2){12.5}}
  \put(95,10){\line(0,1){15}}
  \put(95,25){\line(3,2){12.5}}
  \put(120,10){\line(-3,-2){12.5}}
  \put(120,10){\line(0,1){15}}
  \put(120,25){\line(-3,2){12.5}}
  \put(7.5,17.5){\oval(8,8)[b]}
  \put(7.5,17.5){\oval(8,8)[tr]}
  \put(8.5,21.5){\vector(-1,0){2}}
  \put(40,17.5){\oval(8,8)[b]}
  \put(40,17.5){\oval(8,8)[tr]}
  \put(41,21.5){\vector(-1,0){2}}
  \put(72.5,17.5){\oval(8,8)[b]}
  \put(72.5,17.5){\oval(8,8)[tr]}
  \put(73.5,21.5){\vector(-1,0){2}}
  \put(107.5,17.5){\oval(8,8)[b]}
  \put(107.5,17.5){\oval(8,8)[tr]}
  \put(108.5,21.5){\vector(-1,0){2}}
  \put(0.5,18.5){$\scriptscriptstyle(-1,-1)$}
  \put(9.5,25){$\scriptscriptstyle(0,1)$}
  \put(9.5,9){$\scriptscriptstyle(1,0)$}
  \put(29,25){$\scriptscriptstyle(0,1)$}
  \put(29,9){$\scriptscriptstyle(1,0)$}
  \put(46,25){$\scriptscriptstyle(1,0)$}
  \put(46,9){$\scriptscriptstyle(0,1)$}
  \put(65,25){$\scriptscriptstyle(1,0)$}
  \put(65,9){$\scriptscriptstyle(0,1)$}
  \put(80.5,23){$\scriptscriptstyle(-1,-1)$}
  \put(85.5,12){$\scriptscriptstyle(-1,-1)$}
  \put(96.5,30.2){$\scriptscriptstyle(0,1)$}
  \put(113,30.2){$\scriptscriptstyle(1,0)$}
  \put(96.5,3.7){$\scriptscriptstyle(1,0)$}
  \put(113,3.7){$\scriptscriptstyle(0,1)$}
  \put(110.5,21.5){$\scriptscriptstyle(-1,-1)$}
  \put(23,17){$\scriptstyle\bcs$}
  \put(55,17){$\scriptstyle\bcs$}
  \put(86,17){$\scriptstyle=$}
\end{picture}
\caption{The omnioriented connected sum $\bC P^2\bcs(-S)\bcs\bC
P^2$.}
\end{figure}
\end{exa}

On the other hand, $[\bC P^1]^2$ is also a generator of
$\varOmega^U_4$. It is represented by $(\bC P^1)^2$ with the standard
complex structure, which has second Chern number $4$ and may certainly
be realised over the square.

Our third example shows a related $4$--dimensional situation in which
the connected sum of the quotient polytopes \emph{does\/} support a
suitable orientation.
\begin{exa}\label{cptmincpt}
Let $\overline{\bC P^2}$ denote the quasitoric manifold determined
by the combinatorial data $(\,\overline{\!\varDelta\,}\!(2),-1)$,
whose quotient polytope is the standard $2$-simplex with opposite
orientation. Every vertex has sign $-1$, and we may construct $\bC
P^2\bcs\overline{\bC P^2}$ as an omnioriented quasitoric manifold
over $\varDelta(2)\bcs\varDelta(2)$.  Figure~$2$ illustrates the
procedure diagramatically, in terms of dicharacteristics and
orientations.
\begin{figure}[h]
\begin{picture}(120,35)
  \put(5,5){\line(0,1){30}}
  \put(5,5){\line(5,3){25}}
  \put(5,35){\line(5,-3){25}}
  \put(40,20){\line(5,3){25}}
  \put(40,20){\line(5,-3){25}}
  \put(65,35){\line(0,-1){30}}
  \put(85,5){\line(0,1){30}}
  \put(85,5){\line(1,0){30}}
  \put(85,35){\line(1,0){30}}
  \put(115,5){\line(0,1){30}}
  \put(15,20){\oval(10,10)[b]}
  \put(15,20){\oval(10,10)[tr]}
  \put(15.7,25){\vector(-1,0){2}}
  \put(55,20){\oval(10,10)[b]}
  \put(55,20){\oval(10,10)[tr]}
  \put(55.7,25){\vector(-1,0){2}}
  \put(100,20){\oval(10,10)[b]}
  \put(100,20){\oval(10,10)[tr]}
  \put(100.7,25){\vector(-1,0){2}}
  %
\put(5.3,21.5){$\scriptstyle(-1,-1)$}
  \put(16,29){$\scriptstyle(0,1)$}
  \put(16,9.6){$\scriptstyle(1,0)$}
  \put(47,29){$\scriptstyle(0,1)$}
  \put(47,10){$\scriptstyle(1,0)$}
  \put(65.4,28){$\scriptstyle(-1,-1)$}
  \put(74,10){$\scriptstyle(-1,-1)$}
  \put(98,2){$\scriptstyle(1,0)$}
  \put(98,36){$\scriptstyle(0,1)$}
\put(104,28){$\scriptstyle(-1,-1)$}
  \put(33,19){$\bcs$}
  \put(73,19){$=$}
\end{picture}
\vspace{-3mm} \caption{The omnioriented connected sum $\bC
P^2\bcs\overline{\bC P^2}$.} \label{cp2cs}
\end{figure}

Of course $[\overline{\bC P^2}]=-[\bC P^2]$. So $[\bC
P^2]+[\overline{\bC P^2}]=0$ in $\varOmega_4^U$, and the resulting
manifold bounds by Proposition \ref{scscs}.
\end{exa}

One other observation on $2$--dimensional box sums is also worth
making. Given $k'$-- and $k''$--gons $P'$ and $P''$ in $\bR^2$, it
follows from \eqref{csstats} that
\[
q(P'\bbox P'')=q(P')+q(P'')\sands m(P'\bbox P'')=m(P')+m(P'').
\]
Thus $q(P'\bbox P'')=m(P'\bbox P'')=k'+k''$. So $P'\bbox P''$ is a
$(k'+k'')$--gon, and is combinatorially equivalent to the Minkowski sum
$P'+P''$ whenever $P'$ and $P''$ are in general position.

A situation similar to that of Example \ref{cptboxcpt} arises in
higher dimensions, when we consider the problem of representing
complex cobordism classes by non-singular projective toric varieties. For
any such $V$, the top Chern number coincides with the Euler
characteristic, and is therefore equal to the number of vertices of
the quotient polytope $P$; so $q_-(V)=0$, by \eqref{taraseq}.
Moreover, the Todd genus satisfies $\mathit{Td}(V)=1$.

Omnioriented quasitoric manifolds with $q_-(V)=0$ form an interesting
generalisation of non-singular projective toric varieties, as shown by
Example \ref{mrs}.

\begin{rem}\label{conc1}
Suppose that smooth projective toric varieties $V_1$ and $V_2$ are
of dimension $\geq 4$, and have quotient polytopes $P_1$ and $P_2$
respectively.  Then $c_n(V_1)=q(P_1)$ and $c_n(V_2)=q(P_2)$, yet
$q(P_1\bcs P_2)=q(P_1)+q(P_2)-2$, from \eqref{csstats}. Since
$c_n$ is additive, no omnioriented quasitoric manifold over
$P_1\bcs P_2$ can possibly represent $[V_1]+[V_2]$. This objection
vanishes for $P_1\bbox P_2$, because it enjoys an additional
$2^n-2$ vertices.

The fact that no smooth projective toric variety can represent
$[V_1]+[V_2]$ follows immediately from the Todd genus.
\end{rem}

\begin{exa}\label{mrs}
For any non-negative integers $r$ and $s$ such that $r+s>0$, the
cobordism class $r[\bC P^2]+s[\bC P^1]^2$ is represented by an
omnioriented quasitoric manifold $M(r,s)$. Its quotient polytope
is the iterated box sum
\[
P(r,s)\;=\;\big(\bbox^r\varDelta(2)\big)\bbox\big(\bbox^s I^2\big),
\]
and $q_-(M(r,s))=0$. Applying the Todd genus once more, we deduce that
$M(r,s)$ cannot be cobordant to any smooth toric variety, so long as
$(r,s)\neq(1,0)$ or $(0,1)$.
\end{exa}
%
%
%
%
%
%
%
%
%


\end{document}